\documentclass[a4paper,11pt]{article}
\usepackage{amsmath,amssymb}
\newenvironment{remark}{\vspace{1ex}{\bf Remark.}\rm}{\vspace{1ex}}
\newenvironment{proof}{{\bf Proof.}}{\par\hspace{25em}\rule{1ex}{1ex}\par}
\newtheorem{theorem}{Theorem}[section]
\newtheorem{definition}{Definition}[section]
\newtheorem{corollary}{Corollary}[section]
\newtheorem{lemma}{Lemma}[section]
\newtheorem{proposition}{Proposition}[section]

\title{On special types of minimal and totally geodesic unit vector fields}
\author{Yampolsky A.}
\date{}
\begin{document}
\maketitle

\begin{abstract}
We present a new equation with respect to a unit vector field on
Riemannian manifold $M^n$ such that its solution defines a totally
geodesic submanifold in the unit tangent bundle with Sasaki metric
and apply it to some classes of unit vector fields. We introduce a
class of covariantly normal unit vector fields and prove that
within this class the Hopf vector field is a unique global one
with totally geodesic property. For the wider class of geodesic
unit vector fields on a sphere we give a new necessary and
sufficient condition to generate a totally geodesic submanifold in
$T_1S^n$.


{\it AMS subject class:} Primary 53B20, 53B25; Secondary 53C25.

\end{abstract}

\section*{Introduction}
This paper is organized as follows. In Section 1 we give
definitions of harmonic and minimal unit vector  fields, rough
Hessian and harmonicity tensor for the unit vector field. In
Section 2 we give definition of a totally geodesic unit vector
field and prove the basic Lemma 2.2 which gives a necessary and
sufficient condition for the unit vector field to be totaly
geodesic. The Theorem 2.3 contains a necessary and sufficient
condition on strongly normal unit vector field to be minimal. In
Section 3 we we apply the Lemma 2.2 to the case of a unit sphere
(Lemma 3.1) and describe the geodesic unit vector fields on the
sphere with totally geodesic property (Theorem 3.2). We also
introduce a notion of covariantly normal unit vector field and
prove that within this class the Hopf vector field is a unique one
with totally geodesic property (Theorem 3.1). This theorem is a
revised and simplified version of Theorem 2.1 from \cite{Acta}.
The Section 4 contains an observation that the Hopf vector field
on a unit sphere provides an example of global imbedding of
Sasakian space form into Sasakian manifold as a Sasakian space
form (Theorem 4.1).
\section{Some preliminaries}
\subsection{Sasaki metric}
    Let $(M,g)$ be $n$-dimensional Riemannian manifold  with metric $g$. Denote by
    $\big<\cdot,\cdot\big>$ a scalar product with respect to $g$.
    A natural Riemannian metric on the tangent bundle has been defined by S.~Sasaki \cite{S}.
    We describe it briefly in terms of the {\it connection map}.

    At each point $Q=(q,\xi)\in TM$ the tangent space $T_QTM$ can be split
    into the so-called {\it vertical} and {\it horizontal} parts:
    $$
    T_{Q}TM=\mathcal{H}_{Q}TM \oplus \mathcal{V}_{Q}TM.
    $$
    The vertical part  $\mathcal{V}_{Q}TM$ is tangent to the fiber, while the
    horizontal part is transversal to it.  If $(u^1,\dots,u^n;\xi^1,\dots,\xi^n)$ form
    the natural induced local coordinate system on $TM$, then for $\tilde X \in T_{Q}TM^n$ we
    have
$$
   \tilde X=\tilde X^i \partial /\partial u^i + \tilde X^{n+i} \partial /\partial\xi^i
$$
    with respect to the natural frame $\{ \partial /\partial u^i, \partial /\partial
    \xi^i \}$ on $TM$.

Denote by  $\pi:TM \to M$ the tangent bundle projection map. Then
its differential $\pi_*:T_{ Q}TM \to T_qM $ acts on $\tilde X$ as
$\displaystyle \pi_*\tilde X=\tilde X^i \partial /\partial x^i$
and defines a linear isomorphism between $\mathcal{V}_{ Q}TM$ and
$T_qM$.

The so-called {\it connection map} $K: T_{ Q}TM \to T_qM$ acts on
$\tilde X$ by the rule $\displaystyle K\tilde X=(\tilde
X^{n+i}+\Gamma_{jk}^{i}\xi^j\tilde X^k) \partial /\partial u^i $
and defines a linear isomorphism between $\mathcal{H}_{Q}TM$ and
$T_qM$. The  images $\pi_*\tilde X$ and  $K\tilde X$ are called
{\it horizontal} and {\it vertical } projections of  $\tilde X$,
respectively.  It is easy to see that $\mathcal{V}_{
Q}=\ker\pi_*|_{ Q},\ \mathcal{H}_{Q}=\ker K|_{Q}$.

Let $\tilde X,\tilde Y \in T_{ Q}TM.$ The {\it Sasaki metric} on $TM$ is defined by the following scalar product
    $$
        \big<\big< \tilde X,\tilde Y \big>\big>\big|_Q=
        \big<\pi_* \tilde X, \pi_* \tilde Y\big>\big|_q+\big<K \tilde X,K \tilde Y\big>\big|_q
    $$
    at each point $Q=(q,\xi)$.
    Horizontal and vertical subspaces are mutually orthogonal with respect to Sasaki
    metric.

    The operations inverse to projections are called {\it lifts}. Namely, if
    $X \in T_qM^n$, then
    $\displaystyle X^h=X^i \partial /\partial u^i -\Gamma_{jk}^i\xi^j X^k \partial /\partial \xi^i$
    is in $\mathcal{H}_{Q}TM$ and is called the {\it horizontal lift } of  X, and
    $ \displaystyle X^v=X^i \partial /\partial \xi^i $
    is in $\mathcal{V}_{Q}TM$ and is called the {\it vertical lift } of  $ X$.

    The Sasaki metric can be completely defined by scalar product of
    combinations of lifts of vector fields from $M$ to $TM$ as
    $$
   \big<\big<X^h,Y^h\big>\big>\big|_Q=\big<X,Y\big>\big|_q, \ \
   \big<\big<X^h,Y^v\big>\big>\big|_Q=0,             \ \
    \big<\big<X^v,Y^v\big>\big>\big|_Q=\big<X,Y\big>\big|_q.
   $$
\subsection{Harmonic and minimal unit vector fields}

Suppose, as above, that $u:=(u^1,\dots,u^n)$ are the local coordinates on $M^n$. Denote by
$(u,\xi):=(u^1,\dots,u^n;\xi^1,\dots,\xi^n)$ the natural local coordinates in the tangent bundle $TM^n$. If
$\xi(u)$ is a (unit) vector field on $M^n$, then it defines a mapping
$$
\xi:M^n\to TM^n  \quad\mbox{ or }\quad \xi:M^n\to T_1M^n, \quad \mbox{ if $|\xi|=1$,}
$$
given by $ \xi(u)=(u,\xi(u)).$

For the mappings $f:(M,g)\to (N,h)$ between Riemannian manifolds
the \emph{energy } of $f$ is defined as
$$
E(f):=\frac12 \int_M |d\,f|^2\,d\,Vol_M,
$$
where $|d\,f|$ is a norm of 1-form $d\,f$ in the co-tangent bundle
$T^*M$. Supposing on $T_1M$ the Sasaki metric, the following
definition becomes natural.
\begin{definition}
A unit vector field is called \textit{harmonic}, if it is a critical point of energy functional of mapping
$\xi:M^n\to T_1M^n$.
\end{definition}

Up to an additive constant, the energy functional of the mapping  the  is a total bending of a unit vector field
\cite{Weigm}
$$
B(\xi):=c_n\int_M|\nabla\xi\,|^2\,d\,Vol_M,
$$
where $c_n$ is some normalizing constant and $|\nabla\xi|^2=\sum_{i=1}^n|\nabla_{e_i}\xi|^2$ with respect to
orthonormal frame $e_1,\dots e_n$.

Introduce  a point-wise linear  operator  $A_\xi:T_qM^n\to \xi^\perp_q$, acting as
$$
A_\xi X=-\nabla_X\xi.
$$
In case of integrable distribution $\xi^\perp$, the unit vector
field $\xi$ is called \emph{holonomic} \cite{Am}. In this case the
operator $A_\xi$ is symmetric and is known as Weingarten or a
\emph{shape operator} for each hypersurface of the foliation. In
general, $A_\xi$ is not symmetric, but formally preserves the
Codazzi equation. Namely, a covariant derivative of $A_\xi$ is
defined by
\begin{equation}\label{derA}
-(\nabla_X A_\xi) Y=\nabla_X\nabla_Y\xi-\nabla_{\nabla_XY}\xi.
\end{equation}
Then for the curvature operator of $M^n$ we can write down the Codazzi-type equation
$$
R(X,Y)\xi=(\nabla_Y A_\xi) X-(\nabla_X A_\xi) Y.
$$
From this viewpoint, it is natural to call the operator $A_\xi$ as
\emph{non-holonomic shape} operator.  Remark, that the right hand
side is, up to constant, a \textit{skew symmetric part} of
covariant derivative of $A_\xi$.

Introduce a symmetric tensor field
\begin{equation}\label{Hess}
Hess_\xi(X,Y)=\frac12\big[(\nabla_Y A_\xi) X+(\nabla_X A_\xi) Y\big],
\end{equation}
which is a \textit{symmetric part} of covariant derivative of $A_\xi$. The trace
$$
-\sum_{i=1}^n Hess_\xi(e_i,e_i):=\Delta \xi,
$$
where $e_1,\dots e_n$ is an orthonormal frame, is known as \textit{rough Laplacian} \cite{Besse} of the field
$\xi$. Therefore, one can treat the tensor field \eqref{Hess} as a \textit{rough Hessian} of the field.

With respect to given above notations,  the unit vector field is harmonic if and only if \cite{Weigm}
$$
\Delta\xi=-|\nabla\xi|^2\xi.
$$

Introduce a tensor field
\begin{equation}\label{Hm}
Hm_\xi(X,Y)=\frac12\big[R(\xi,A_\xi X)Y+R(\xi,A_\xi Y)X\big],
\end{equation}
which is a symmetric part of tensor field $R(\xi,A_\xi X)Y$. The
trace
$$
trace\, Hm_\xi:=\sum_{i=1}^n Hm_\xi(e_i,e_i)
$$
is responsible for harmonicity  of mapping $\xi:M^n\to T_1M^n$ in
terms of general notion of harmonic maps \cite{E-L}. Precisely, a
\emph{harmonic} unit vector field $\xi$ defines a \emph{harmonic
mapping} $\xi:M^n\to T_1M^n$ if and only if  \cite{GM}
$$
trace\, Hm_\xi=0.
$$
From this viewpoint, it is natural to call the tensor field \eqref{Hm} as \textit{harmonicity tensor} of the
field $\xi$.

Consider now the image $\xi(M^n)\subset T_1M^n$ with a pull-back Sasaki metric.
\begin{definition} A unit vector field  $\xi$  on Riemannian manifold $M^n$ is called minimal if the
image of (local) imbedding $\xi:M^n\to T_1M^n$ is minimal submanifold in the unit tangent bundle $T_1M^n$ with
Sasaki metric.
\end{definition}
A number of results on minimal unit vector fields one can find in
\cite{BX-V2,BX-V1,Bx-Vh4,B-Ch-N, GM-LF2, GM-GD-Vh, GM-LF, G-Z,
GD-V1,GD-Vh2, Ped, TS-Vh2, TS-Vh3,TS-Vh1}. In \cite{Ym1}, the
author has found explicitly the second fundamental form of
$\xi(M^n)$ and presented some examples of unit vector fields of
\emph{constant mean curvature}.

\section{Totally geodesic unit vector fields}

\begin{definition} A unit vector field  $\xi$  on Riemannian manifold $M^n$ is called totally geodesic if the
image of (local) imbedding $\xi:M^n\to T_1M^n$ is totally geodesic submanifold in the unit tangent bundle
$T_1M^n$ with Sasaki metric.
\end{definition}

Using the explicit expression for the second fundamental form
\cite{Ym1}, the author gave a full description of the totally
geodesic (local) unit vector fields on 2-dimensional Riemannian
manifold.
\begin{theorem} \cite{Ym5}
Let $(M^2,g)$ be a Riemannian manifold with sign-preserving
Gaussian curvature $K$. Then $M$ admits a totally geodesic unit
vector field $\xi$ if and only if there is a local parametrization
of $M$ with respect to which the metric $g$ is of the form
$$
ds^2=du^2+\sin^2 \alpha(u) \, dv^2,
$$
where $\alpha(u)$ solves the differential equation
$\quad\displaystyle \frac{d\alpha}{du}=1-\frac{a+1}{\cos\alpha}$.
The corresponding local unit vector field  $\xi$ is of the form
$$
\xi=\cos (av+\omega_0)\,\partial_u+\frac{\sin
(av+\omega_0)}{\sin\alpha(u)}\,\partial_v,
$$
where $a,\omega_0=const$.
\end{theorem}
For the case of \emph{flat} Riemannian 2-manifold, {the totally
geodesic unit vector field is either parallel or moves helically
along a pencil of parallel straight lines on a plane with a
constant angle speed } \cite{Ym2} . It is easy to see that the
following corollary is true.
\begin{corollary}
Integral trajectories of a totally geodesic (local) unit vector
field on the non-flat Riemannian manifold $M^2$ are locally
conformally equivalent to the integral trajectories of totally
geodesic unit vector field on a plane. Moreover, with respect to
Cartesian coordinates $(x,y)$ on the plane, these integral
trajectories are
$$
\begin{array}{ll}
x=c& \mbox{for a=0,}\\[1ex]
y(x)=-\frac{1}{a}\ln |\sin(ax)|+c &\mbox{for $a\ne 0$,}
\end{array}
$$
where $c$ is a parameter.
\end{corollary}

In what follows,  we present a new differential equation with
respect to a unit vector field such that its solution  generates a
totally geodesic submanifold in $T_1M^n$.

In terms of horizontal and vertical lifts of vector fields from
the base to its tangent bundle, the differential of mapping
$\xi:M^n\to TM^n$ is acting as
\begin{equation}\label{Eq0}
\xi_*X=X^h+(\nabla_X \xi)^v=X^h-(A_\xi X)^v,
\end{equation}
where $\nabla$ means Levi-Civita connection on $M^n$ and the lifts
are considered to points of $\xi(M^n)$.

It is well known that if $\xi$ is a \emph{unit} vector field on
$M^n$, then the vertical lift $\xi^v$ is a \textit{unit normal}
vector field on a hypersurface $T_1M^n\subset TM^n$. Since $\xi$
is of unit length, $\xi_* X\perp \xi^v$ and hence in this case $
\xi_*:TM^n\to T(T_1M^n)$.

Denote by $A_\xi^t:\xi_q^\perp\to T_qM^n$ a formal adjoint
operator
$$
\big<A_\xi X,Y\big>_q=\big<X, A^t_\xi Y\big>_q.
$$
Denote by $\xi^\perp$ a distribution on $M^n$ with $\xi$ as its
normal unit vector field. Then for each vector field $N\in
\xi^\perp$, the vector field
\begin{equation}\label{normal}
\tilde N=(A_\xi^t N)^h+N^v
\end{equation}
is normal to $\xi(M^n)$. Thus, \eqref{normal} presents the normal
distribution on $\xi(M^n)$.
\begin{lemma}\label{main_Lemma}
Let $M^n$ be Riemannian manifold and $T_1M^n$ its unit tangent
bundle with Sasaki metric. Let $\xi$ a smooth (local) unit vector
field on $M^n$. The second fundamental form $\tilde\Omega_{\tilde
N}$ of $\xi(M^n)\subset T_1M^n$ with respect to the normal vector
field \eqref{normal} is of the form
\begin{equation}\label{Form}
\tilde\Omega_{\tilde N}(\xi_* X,\xi_* Y)
=-\big<Hess_\xi(X,Y)+A_\xi Hm_\xi(X,Y),N\big>,
\end{equation}
where $X$ and $Y$ are arbitrary vector fields  on $M^n$.
\end{lemma}
\begin{proof}
By definition, we have
$$ \tilde\Omega_{\tilde N}(\xi_* X,\xi_* Y)=\big<\big<\tilde \nabla_{\xi_*X}\,\xi_*Y,\tilde
N\big>\big>_{(q,\xi(q))},
$$
where $\tilde \nabla$ is the Levi-Civita connection of Sasaki metric on $TM^n$. To calculate $\tilde
\nabla_{\xi_*X}\,\xi_*Y$, we  can use the  formulas \cite{Kow}
\begin{equation}\label{Kow}
\begin{array}{ll}
        \tilde\nabla_{X^h}Y^h  =  (\nabla_XY)^h - \frac{1}{2}(R(X,Y)\xi)^v,
         &\tilde\nabla_{X^v}Y^h  =  \frac{1}{2}(R(\xi,X)Y)^h, \\[1ex]
        \tilde\nabla_{X^h}Y^v  =  (\nabla_XY)^v \ + \frac{1}{2}(R(\xi,Y)X)^h, &
        \tilde\nabla_{X^v}Y^v  =  0.
\end{array}
\end{equation}
A direct calculation yields
\begin{multline*}\label{der}
$$
\tilde \nabla_{\xi_*X}\,\xi_*Y=\Big( \nabla_XY+\frac12R(\xi,\nabla_X\xi)Y+\frac12R(\xi,\nabla_Y\xi)X\Big)^h+\\
\Big(\nabla_X\nabla_Y\xi-\frac12R(X,Y)\xi\Big)^v.
$$
\end{multline*}
The derivative above is not tangent to $\xi(M^n)$. It contains a
projection on "external" normal vector field, i.e. on $\xi^v$
which is a unit normal of $T_1M^n$ inside $TM^n$. To correct the
situation, we should subtract this projection, namely
$-\big<\nabla_X\xi,\nabla_Y\xi\big>\xi$, from the vertical part of
the derivative.

 Therefore, we have
\begin{multline*}
\tilde\Omega_{\tilde N}(\xi_* X,\xi_* Y)=\big<\nabla_X\nabla_Y\xi+\big<\nabla_X\xi,\nabla_Y\xi\big>\xi-\frac12R(X,Y)\xi,N\big>+\\
\big<\nabla_XY+\frac12R(\xi,\nabla_X\xi)Y+\frac12R(\xi,\nabla_Y\xi)X,A_\xi^tN\big>
\end{multline*}
or, equivalently,
\begin{multline*}
\tilde\Omega_{\tilde N}(\xi_* X,\xi_* Y)=\big<\nabla_X\nabla_Y\xi+\big<\nabla_X\xi,\nabla_Y\xi\big>\xi-\frac12R(X,Y)\xi+\\
A_\xi\big(\nabla_XY+\frac12R(\xi,\nabla_X\xi)Y+\frac12R(\xi,\nabla_Y\xi)X\big),N\big>.
\end{multline*}
Taking into account \eqref{derA}, \eqref{Hess}, \eqref{Hm} and
\eqref{normal}, and also
$R(X,Y)\xi=\nabla_X\nabla_Y\xi-\nabla_Y\nabla_X\xi-\nabla_{[X,Y]}\xi$,
we can write
$$
\tilde\Omega_{\tilde N}(\xi_* X,\xi_* Y)=-\big<Hess_\xi(X,Y)+A_\xi
Hm_\xi(X,Y),N\big>
$$
which completes the proof.
\end{proof}

\begin{lemma}
Let $M^n$ be Riemannian manifold and $T_1M^n$ its unit tangent
bundle with Sasaki metric. Let $\xi$ be a smooth (local) unit
vector field on $M^n$. The vector field $\xi$ generates a  totally
geodesic submanifold $\xi(M^n)\subset T_1M^n$ if and only if $\xi$
satisfies
\begin{equation}\label{MainEq}
Hess_\xi(X,Y)+A_\xi Hm_\xi(X,Y)-\big<A_\xi X,A_\xi Y\big>\,\xi=0
\end{equation}
for all (local) vector fields $X,Y$ on $M^n$.
\end{lemma}
\begin{proof}
Taking into account \eqref{Form}, the condition on $\xi$ to be
totally geodesic takes the form
$$
-Hess_\xi(X,Y)-A_\xi Hm_\xi(X,Y)=\lambda\,\,\xi.
$$
Multiplying the equation above by $\xi$, we can find easily
$\lambda=-\big<A_\xi X,A_\xi Y\big>$.
\end{proof}
Follow \cite{GD-V1}, we call a unit vector field $\xi$
\emph{strongly normal} if
$$
\big<(\nabla_XA_\xi)Y,Z\big>=0
$$
for all $X,Y,Z\in\xi^\perp$. In other words, $
(\nabla_XA_\xi)Y=\lambda\xi$ for all $X,Y\in \xi^\perp$. It is
easy to find the function $\lambda$.  Indeed, we have
$$
\lambda=\big<(\nabla_XA_\xi)Y,\xi\big>=\big<\nabla_{\nabla_XY}\xi-\nabla_X\nabla_Y\xi,\xi\big>=-\big<\nabla_X\nabla_Y\xi,\xi\big>=
\big<\nabla_X\xi,\nabla_Y\xi\big>.
$$
Thus, the strongly normal unit vector  field can be characterized
by the equation
\begin{equation}\label{StrN}
(\nabla_XA_\xi)Y=\big<A_\xi X,A_\xi Y\big>\,\xi
\end{equation}
for all $X,Y\in \xi^\perp$.

The strong normality condition highly simplifies the second
fundamental form of $\xi(M^n)\subset T_1M^n$. An orthonormal frame
$e_1,e_2,\dots,e_n$ is called {\it adapted} to the field $\xi$  if
$e_1=\xi$ and $e_2,\dots,e_n\in \xi^\perp$.
\begin{proposition}\label{2.1} Let $\xi$ be a unit strongly normal vector
field on Riemannian manifold $M^n$. With respect to the adapted
frame, the matrical components of the second fundamental form of
$\xi(M^n)\subset T_1(M^n)$ simultaneously take the form
$$
\tilde\Omega_{\tilde N}=
\begin{pmatrix}
 * & * & \dots & * \\
 * & 0 &\dots  & 0 \\
 \vdots &\vdots  & & \vdots\\
  * & 0 &\dots &0
\end{pmatrix}.
$$
\end{proposition}
\begin{proof}
Set $N_\sigma=e_\sigma\quad (\sigma=2,\dots,n)$. The condition
\eqref{StrN} implies
$$
R(X,Y)\xi=0,\quad Hess_\xi(X,Y)=\big< A_\xi X,A_\xi Y\big>
\xi,\quad Hm_\xi(X,Y)\,\sim\, \xi
$$
for all $X,Y\in \xi^\perp$. Therefore, with respect to the adapted
frame $$ \tilde \Omega_{\sigma}(\xi_* e_\alpha,\xi_*
e_\beta)=0\quad (\alpha,\beta=2,\dots, n)
$$
for all $\sigma=2,\dots,n$.
\end{proof}
The following assertion is a natural corollary of the Proposition
\ref{2.1} .
\begin{theorem} Let $\xi$ be a unit strongly normal vector field. Denote by $k$ the geodesic curvature of
its integral trajectories and  by $\nu$ the principal normal unit
vector field of the trajectories. The field  $\xi$ is minimal if
and only if
$$
k[\xi,\nu]+\xi(k)\nu -k A_\xi R(\nu,\xi)\xi+k^2 \xi=0
$$
where $[\xi,\nu]=\nabla_\xi\nu-\nabla_\nu\xi$.
\end{theorem}
\begin{proof}
Indeed,
$$
\Tilde
\Omega_{\sigma}(\xi_*e_1,\xi_*e_1)=-\big<Hess_\xi(\xi,\xi)+A_\xi
Hm_\xi(\xi,\xi),e_\sigma\big>
$$
Denote by $\nu$ a vector field of the principal normals of
$\xi$-integral trajectories and by $k$ their geodesic curvature
function. Then
$$
\begin{array}{l}
Hess_\xi(\xi,\xi)=\nabla_{\nabla_\xi\xi}-\nabla_\xi\nabla_\xi\xi=k\nabla_\nu\xi-\nabla_\xi(k\nu)=k[\nu,\xi]-\xi(k)\nu,
\\[2ex]
Hm_\xi(\xi,\xi)=-R(\xi,\nabla_\xi\xi)\xi=-kR(\xi,\nu)\xi
\end{array}
$$
and we get
$$
\Tilde \Omega_{\sigma}(\xi_*e_1,\xi_*
e_1)=\big<k[\xi,\nu]+\xi(k)\nu -k A_\xi
R(\nu,\xi)\xi,e_\sigma\big>.
$$
Finally, to be minimal, the field $\xi$ should satisfy
$$
k[\xi,\nu]+\xi(k)\nu -k A_\xi R(\nu,\xi)\xi=\lambda\,\xi.
$$
Multiplying by $\xi$, we get
$$
\lambda=k\big<[\xi,\nu],\xi\big>=k\big<\nabla_\xi\nu,\xi\big>=-k^2,
$$
which completes the proof.
\end{proof}
Thus, we get the following.
\begin{corollary}\cite{GD-V1}
Every unit strongly normal geodesic vector field is minimal.
\end{corollary}
Most of examples of minimal unit vector fields in \cite{GD-V1} are
based on this Corollary.

\section{The case of a unit sphere}

If the manifold is a unit sphere $S^{n+1}$, the equation
\eqref{MainEq} can be essentially simplified.
\begin{lemma}\label{2}
A unit  (local) vector field $\xi$ on a unit sphere  $S^{n+1}$ generates a totally geodesic submanifold
$\xi(S^{n+1})\subset T_1S^{n+1}$ if and only if $\xi$ satisfies
\begin{equation}\label{New_Eq}
\begin{array}{r}
(\nabla_X A_\xi)Y=\dfrac12\Big[(\mathcal{L}_\xi \,g)(X,Y)\,A_\xi\xi+\big<\xi,X\big>\,(A^2_\xi Y+Y)+\\
\big<\xi,Y\big>\,(A^2_\xi X-X)\Big] + \big<A_\xi X,A_\xi Y\big>\,\xi,
\end{array}
\end{equation}
where $(\mathcal{L}_\xi\, g)(X,Y)=\big<\nabla_X\xi,Y\big>+\big<X,\nabla_Y\xi\big>$ is a Lie derivative of metric
tensor in a direction of $\xi$.
\end{lemma}
\begin{proof}
Indeed, on a unit sphere
$$
(\nabla_Y A_\xi) X-(\nabla_X A_\xi) Y=R(X,Y)\xi=\big<\xi,Y\big>X-\big<\xi,X\big>Y.
$$
Hence,
$$
Hess_\xi(X,Y)=(\nabla_X A_\xi) Y+\frac12[\big<\xi,Y\big>X-\big<\xi,X\big>Y].
$$
For $Hm_\xi(X,Y)$ we have
\begin{multline*}
Hm_\xi(X,Y)=\frac12\Big[\big<\nabla_X\xi,Y\big>\xi-\big<\xi,Y\big>\nabla_X\xi+
\big<\nabla_Y\xi,X\big>\xi-\big<\xi,X\big>\nabla_Y\xi\Big]=\\
\frac12(\mathcal{L}_\xi\,g)(X,Y)\,\xi +\frac12\Big[\big<\xi,Y\big>A_\xi X+\big<\xi,X\big>A_\xi Y \Big].
\end{multline*}
Finally, we find
\begin{multline*}
(\nabla_X A_\xi)Y=\\
\frac12\Big[(\mathcal{L}_\xi\,g)(X,Y)\,A_\xi\xi+\big<\xi,X\big>\,(A^2_\xi
Y+Y)+\big<\xi,Y\big>\,(A^2_\xi X-X)\Big] +\big<A_\xi X,A_\xi
Y\big>\,\xi.
\end{multline*}

\end{proof}

Remind that the operator $A_\xi$ is symmetric  if and only if the
field $\xi$ is holonomic, and is skew-symmetric if and only if the
field $\xi$ is a Killing vector field. Both types of these fields
can be included into a class of \emph{covariantly normal} unit
vector fields.
\begin{definition} A regular unit vector field on Riemannian manifold is said to be covariantly normal if the
operator $A_\xi:TM\to \xi^\perp$ defined by $A_\xi X=-\nabla_X\xi$
satisfies the normality condition
$$
A_\xi^tA_\xi=A_\xi A_\xi^t
$$
with respect to some orthonormal frame.
\end{definition}

The integral trajectories of holonomic and Killing unit vector
fields are always geodesic. Every covariantly normal unit vector
field possesses this property.
\begin{lemma} Integral trajectories of a covariantly normal unit vector field are geodesic lines.
\end{lemma}
\begin{proof}
Suppose $\xi$ is a unit covariantly normal vector field on a Riemannian manifold $M^{n+1}$. Find a unit vector
field $\nu_1$ such that
$$
\nabla_\xi\xi=-k\nu_1.
$$
Geometrically, the function $k$ is a geodesic curvature of the
integral trajectory of the field $\xi$.

Complete up the pair $(\xi,\nu_1)$ to the orthonormal frame $(\xi,\nu_1,\dots \nu_{n})$. Then we can set
$$
\nabla_\xi\xi=-k\nu_1,\quad \nabla_{\nu_\alpha}\xi=-a_\alpha^\beta\nu_\beta,
$$
where $\alpha,\beta=1,\dots,n$. With respect to the frame $(\xi,\nu_1,\dots \nu_{n})$ the matrix $A_\xi$ takes
the form
$$-A_\xi=
\left(%
\begin{array}{cllcl}
  0 & k & 0 & \dots & 0 \\
  0 & a^1_1 & a^1_2 & \dots & a^1_{n} \\
  \vdots & \vdots & \vdots & \vdots & \vdots \\
  0 & a^{n}_1 & a^{n}_2 & \dots & a^{n}_{n} \\
\end{array}%
\right)
$$
and as a consequence
$$-A^t_\xi=
\left(%
\begin{array}{cllcl}
  0 & 0 & 0 & \dots & 0 \\
  k & a^1_1 & a_1^2 & \dots & a_1^{n} \\
  \vdots & \vdots & \vdots & \vdots & \vdots \\
  0 & a_{n}^1 & a_{n}^2 & \dots & a_{n}^{n} \\
\end{array}%
\right).
$$
Therefore,
$$
A_\xi A_\xi^t=\left(
\begin{array}{cccc}
  k^2 & ka_1^1 & \dots & ka^n_1 \\
  ka_1^1 &* &  \dots & * \\
  \vdots &  \vdots & \vdots & \vdots \\
 ka^n_1 & * &  \dots & * \\
\end{array}%
\right), \quad A_\xi^t A_\xi=\left(
\begin{array}{cccc}
  0 & 0 & \dots & 0 \\
  0 &* &  \dots & * \\
  \vdots &  \vdots & \vdots & \vdots \\
 0 & * &  \dots & * \\
\end{array}%
\right)
$$
and we conclude $k=0$.

\end{proof}

Now we can easily prove the following theorem.
\begin{theorem}\label{Hopf}
Let $\xi$ be a global covariantly normal unit vector field on a
unit sphere $S^{n+1}$. Then $\xi$ is a totally geodesic if and
only if $n=2m$ and $\xi$ is a Hopf vector field.
\end{theorem}
\begin{proof}
Suppose $\xi$ is covariantly normal and totally geodesic. Then $$A_\xi\xi=-\nabla_\xi\xi=0$$ by Lemma \ref{2}
and the equation \eqref{New_Eq} takes the form
\begin{equation}\label{New_Eq1}
(\nabla_X A_\xi)Y=\dfrac12\Big[\big<\xi,X\big>\,(A^2_\xi Y+Y)+\big<\xi,Y\big>\,(A^2_\xi X-X)\Big] + \big<A_\xi
X,A_\xi Y\big>\,\xi.
\end{equation}
Setting $X=Y=\xi$ we get an identity. Set $Y=\xi$ and take arbitrary unit $ X\perp \xi$. Then we get
$$
2(\nabla_X A_\xi)\xi+X=A_\xi^2X.
$$
On the other hand, directly
$$
(\nabla_X A_\xi)\xi=-(\nabla_X\nabla_\xi\xi-\nabla_{\nabla_X\xi}\xi)=A_\xi^2X.
$$
Hence,
$$
A_\xi^2\big|_{\xi^\perp}=-E.
$$
Therefore, $n=2m$.  Since $A_\xi$ is real normal linear operator, there exists an orthonormal frame such that
$$
A_\xi=
\left(%
\begin{array}{cccc}
0&  &  &  \\
 & \begin{array}{rr}
    0 & 1 \\
    -1 & 0 \\
  \end{array} &  &   \\
&   & \ddots &    \\
&   &  &   \begin{array}{rr}
     0 & 1 \\
     -1 & 0 \\
   \end{array} \\
\end{array}%
\right)
$$
with zero all other entries. Therefore, $A_\xi+A_\xi^t=0$ and $
\xi$ is a Killing vector field. Since $\xi$ is supposed global,
$\xi$ is a Hopf vector field.

Finally, if we take  $X,Y\perp\xi$, we get the equation
$$
(\nabla_X A_\xi)Y=\big<A_\xi X,A_\xi Y\big>\xi.
$$
But for a Killing vector field $\xi$ we have \cite{GD-V1}
$$
(\nabla_X A_\xi)Y=R(\xi,X)Y=\big<X,Y\big>\xi.
$$
Since $\xi$ is a Hopf vector field,
$
\big<A_\xi X,A_\xi Y\big>=\big<X,Y\big>.
$
So, in this case we have an identity.

If we suppose now that $\xi$ is a Hopf vector field on a unit sphere, then $\xi$ is covariantly normal as a
Killing vector field and totally geodesic  \cite{Acta} as a characteristic vector field of a standard contact
metric structure on $S^{2m+1}$.
\end{proof}

\begin{remark}
Theorem \ref{Hopf} is a correct and simplified version of Theorem 2.1 \cite{Acta}, where the normality of the
operator $A_\xi$ was implicitly used in a proof.
\end{remark}

In the case of a \emph{weaker condition} on the field $\xi$ to be
only a \emph{geodesic} one, the result is not so definite. We
begin with some preparations.

The almost complex structure on $TM^n$ is defined by
$$
JX^h=X^v,\quad JX^v=-X^h
$$
for all vector field $X$ on $M^n$. Thus, $TM^n$ with Sasaki metric
is an almost K\"ahlerian manifold. It is K\"ahlerian if and only
if $M^n$ is flat \cite{Dmb}.

The unit tangent bundle $T_1M^n$ is a hypersurface in $TM^n$ with
a unit normal vector $\xi^v$ at each point $(q,\xi)\in T_1M^n$.
Define a unit vector field $\bar \xi$, a 1-form $\bar\eta$ and a
$(1,1)$ tensor field $\bar\varphi$ on $T_1M^n$ by
$$
\bar\xi=-J\xi^v=\xi^h, \quad JX=\bar\varphi X+\bar\eta(X)\xi^v.
$$
The triple $(\bar\xi,\bar\eta,\bar\varphi)$ form a standard almost
contact structure on $T_1M^n$ with Sasaki metric $g_S$. This
structure is not almost contact \emph{metric} one. By taking
$$
\tilde\xi=2\bar\xi=2\xi^h, \quad \tilde\eta=\frac12\bar\eta,\quad
\tilde\varphi=\bar\varphi,\quad g_{cm}=\frac14g_S
$$
at each point $(q,\xi)\in T_1M^n$, we get the \emph{almost contact
metric structure} $(\tilde\xi,\tilde\eta,\tilde\varphi)$ on
$(T_1M^n,g_{cm})$.

In a case of a general almost contact metric manifold $(\tilde
M,\tilde \xi,\tilde\eta,\tilde\varphi,\tilde g)$ the following
definition is known \cite{Blair}.
\begin{definition}
A submanifold $N$ of a contact metric manifold $(\tilde M,\tilde
\xi,\tilde\eta,\tilde\varphi,\tilde g)$ is called invariant if
$\tilde\varphi(T_pN)\subset T_pN$ and anti-invariant if
$\tilde\varphi(T_pN)\subset (T_pN)^\perp$ for every $p\in N$.
\end{definition}

If $N$ is the invariant submanifold, then the characteristic
vector field $\tilde \xi$ is \textbf{tangent} to $N$ at each of
its points.

After all mentioned above, the following definition is natural
\cite{Binh}.
\begin{definition}
A unit vector field $\xi$ on a Riemannian manifold $(M^n,g)$ is
called invariant (anti-invariant) is the submanifold
$\xi(M^n)\subset (T_1M^n,g_{cm})$ is invariant (anti-invariant).
\end{definition}

It is easy to see from \eqref{Eq0} that the \emph{invariant} unit
vector field is always a geodesic one, i.e. its integral
trajectories are geodesic lines.

Binh T.Q., Boeckx E. and Vanhecke L. have considered this kind of
unit vector fields \cite{Binh} and proved the following
proposition.
\begin{proposition}\label{Binh}
A unit vector field $\xi$ on $(M^n, g)$ is invariant if and only
if $(\tilde\xi=\xi\,,\, \tilde \eta=\big<\cdot,\xi\big>_g\,,\,
\tilde\varphi=A_\xi)$ is an almost contact structure on $M^n$. In
particular, $\xi$ is a geodesic vector field on $M^n$ and
$n=2m+1$.
\end{proposition}

Now we can formulate the result.
\begin{theorem}\label{Tm}
A unit geodesic vector field $\xi$ on $S^{n+1}$ is totally geodesic if and only if $n=2m$ and $\xi$ is a
strongly normal invariant unit vector field.
\end{theorem}
\begin{proof}
Suppose  $\xi$ is a geodesic  and totally geodesic unit vector field. Then $A_\xi\xi=0$ and the equation
\eqref{New_Eq} takes the form \eqref{New_Eq1}. Follow the proof of the Theorem \ref{Hopf}, we come to the
following conditions on the field $\xi$:
\begin{equation}\label{Tg}
A_\xi^2X=-X,\quad (\nabla_X A_\xi)Y= \big<A_\xi X,A_\xi Y\big>\,\xi
\end{equation}
for all $X,Y\in \xi^\perp$. From $\eqref{Tg}_1$ we conclude that $n=2m$. Comparing $\eqref{Tg}_2$ with
\eqref{StrN}, we see that $\xi$ is a strongly normal vector field.

Consider now a $(1,1)$ tensor field $\varphi=A_\xi=-\nabla\xi$ and a 1-form $\eta=\big<\,\cdot\,,\,\xi\,\big>$.
Taking into account $\eqref{Tg}_1$ and $A_\xi\xi=0$, we see that
$$
\varphi^2X=-X+\eta(X)\xi,\quad \varphi\xi=0,\quad\eta(\varphi X)=0,\quad \eta(X)=1
$$
for any vector field $X$ on the sphere.  Therefore, the triple
$$
\tilde \varphi=A\xi,\quad \tilde\xi=\xi,\quad \tilde\eta=\big<\,\cdot\,,\,\xi\,\big>
$$
form an \emph{almost contact structure} with the field $\xi$ as a characteristic vector field of this structure.
By Proposition \ref{Binh}, the field $\xi$ is invariant.

Conversely, suppose $\xi$ is strongly normal and invariant on
$S^{n+1}$.  Then, by Proposition \ref{Binh}, $\xi$ is geodesic and
$n=2m$. The rest of the proof is a direct checking of the formula
\eqref{New_Eq1}.
\end{proof}
\section{A remarkable property of the Hopf vector field}

It is well-known that for a unit sphere $S^n$ the standard contact
metric structure on $T_1S^n$ is a Sasakian one. If $\xi$ is a Hopf
unit vector field on $S^{2m+1}$, then $\xi$ is  a characteristic
vector field of a standard contact metric structure on the unit
sphere $S^{2m+1}$. By Proposition \ref{Binh}, the submanifold
$\xi(S^{2m+1})$ is the invariant submanifold in $T_1S^{2m+1}$.
Therefore, $\xi(S^{2m+1})$ \emph{is also Sasakian} with respect to
the induced structure \cite{Y-K}.  Since the Hopf vector field is
strongly normal, by Theorem \ref{Tm}, the submanifold
$\xi(S^{2m+1})$ is totally geodesic. The sectional curvature of
the submanifold $\xi(S^{2m+1})$ was found in \cite{Acta} and
implies a remarkable corollary.
\begin{theorem}
Let $\xi$ be a Hopf vector field on the unit sphere $S^{2m+1}$. With respect to the induced structure, the
manifold $\xi(S^{2m+1})$ is a Sasakian space form of $\varphi$- curvature 5/4.
\end{theorem}
In other words, the Hopf vector field provides an example of
embedding of a \emph{Sasakian space form} of $\varphi$-curvature 1
into \emph{Sasakian manifold} such that the image is contact,
totally geodesic \emph{Sasakian space form} of $\varphi$-curvature
5/4 with respect to the induced structure.

\vspace{1cm}

\noindent
Department of Geometry,\\
Faculty of Mechanics and Mathematics,\\
Kharkiv National University,\\
Svobody Sq. 4,\\
 61077, Kharkiv,\\
Ukraine.\\
e-mail: yamp@univer.kharkov.ua
\end{document}